\documentclass[12pt]{amsart}
\numberwithin{equation}{section}
\newtheorem{theorem}{Theorem}
\numberwithin{theorem}{section}

\numberwithin{lemma}{section}
\newtheorem{definition}{Definition}
\numberwithin{definition}{section}

\numberwithin{remark}{section}

\numberwithin{example}{section}
\newtheorem{corollary}{Corollary}
\numberwithin{corollary}{section}

\numberwithin{proposition}{section}
\def\b{\begin{equation}}
\def\e{\end{equation}}
\newcommand{\ignore}[1]{}
\normalsize
\date {April 3, 2005}
\thanks{AMS Subject Classifications: 35K65,35R05, 22E30, 26D10 }
\keywords{nonlinear parabolic equations, positive solutions,
 Hardy-type inequality}
\begin{document}
\pagenumbering{arabic} \pagenumbering{arabic}\setcounter{page}{1}
\tracingpages 1
\title{The Hardy inequality and Nonlinear parabolic equations on Carnot groups}
\author{Ismail Kombe}
\address{Mathematics Department\\ Dawson-Loeffler Science
\&Mathematics Bldg\\
Oklahoma City University \\
2501 N. Blackwelder, Oklahoma City, OK 73106-1493}
\email{ikombe@okcu.edu}
\begin{abstract}
In this paper we shall investigate the nonexistence of positive
solutions for the following nonlinear parabolic partial
differential equation:\[
\begin{cases}
\frac{\partial u}{\partial t}= \Delta_{\mathbb{G},p}u+V(x)u^{p-1}
& \text{in}\quad \Omega \times (0, T ), \quad 1<p<2 ,\\
u(x,0)=u_{0}(x)\geq 0 & \text{in} \quad\Omega, \\
u(x,t)=0 & \text{on}\quad
\partial\Omega\times (0, T)
\end{cases}
\]
where $ \Delta_{\mathbb{G},p}$ is the $p$-sub-Laplacian on Carnot
group $ \mathbb{G}$  and $V\in L_{\text{loc}}^1(\Omega)$.
\end{abstract}
\maketitle
\section{Introduction}

In this paper we study the nonexistence of positive solutions  for
the  following nonlinear degenerate parabolic equation:
\begin{equation}
\begin{cases}
\frac{\partial u}{\partial t}=\Delta_{\mathbb{G},p}\,u+V(x)u^{p-1}
& \text{in}\quad \Omega \times (0, T ),\\
u(x,0)=u_{0}(x)\geq 0 & \text{in} \quad\Omega, \\
u(x,t)=0 & \text{on}\quad
\partial\Omega\times (0, T).
\end{cases}
\end{equation}
Here $\Omega$ is a bounded domain with smooth boundary in
$\mathbb{G}$, $V\in L_{\text{loc}}^1(\Omega)$ and $1<p<2$. The
nonlinear operator
\begin{equation}
\Delta_{\mathbb{G},p}\,u=\nabla_{ \mathbb{G}}\cdot (|\nabla_{
\mathbb{G}}u|^{p-2}\nabla_{ \mathbb{G}}u)
\end{equation}
is the sub-$p$-Laplacian on Carnot group $\mathbb{G}$.

 It is well known that Euclidean space $
\mathbb{R}^n$ with its usual Abelian group structure is a trivial
Carnot group and the problem (1.1) has been studied by Goldstein
and Kombe \cite{20} on $ \mathbb{R}^n$. In this case we should
also mention that the problem (1.1) with the potential
$V(x)=\frac{c}{|x|^p}$ has been studied in great detail by Garcia
and Peral \cite{16}, and by Aguilar and Peral \cite{1}.  Recently,
Goldstein and Kombe \cite{21} extended the results in \cite {20}
to the Heisenberg group $ \mathbb{H}^n$ which is the simplest
non-commutative Carnot group.

It turns out that nonexistence of positive solutions that kind of
problems largely depends on the size of bottom of the normalized
$p$-energy form
\begin{equation}
\sigma_{\inf}^p(V):=\inf_{0 \not\equiv \phi \in C_0^{\infty}
(\Omega)}\frac{\int_\Omega|\nabla_{
\mathbb{G}}\phi|^pdx-\int_\Omega V|\phi|^p dx}{\int_\Omega
|\phi|^p dx}
\end{equation}
and the range of $p$. This was discovered by Cabr\'e and Martel
\cite {8} for $p=2$ in Euclidean space $ \mathbb{R}^n$ which was a
nice generalization of Baras and Goldstein  result \cite {4}. The
observation of Cabr\'e and Martel \cite {8} used to get additional
results on nonexistence of positive solutions for wide class of
linear and nonlinear parabolic problems on Euclidean space $
\mathbb{R}^n$ and the Heisenberg group $ \mathbb{H}^n$ ( see \cite
{18}, \cite {19},  \cite {20}, \cite {21}, \cite {23}, \cite {26},
and \cite {27}).

The purpose of this paper is to investigate nonexistence of
positive solutions of the problem(1.1) on more general Carnot
group $ \mathbb{G}$ and  our results recover the previous results
\cite{20}, \cite {21}. We considered
 singular  potentials $V(x)=c/|x|^p$ and
$V(z,l)=c|z|^p/(|z|^4+l^2)^{p/2}$ in \cite {20}  and \cite {21},
respectively. These potentials appear in the Hardy's inequality on
Euclidean space  $ \mathbb{R}^n$ and the Heisenberg group
$\mathbb{H}^n$, respectively. It is natural to investigate such a
singular potential for the problem (1.1) on more general Carnot
groups and we shall first prove a Hardy type  inequality on
polarizable Carnot groups.

The plan of the paper is as follows. In Section 2 we recall the
basic properties of Carnot group $\mathbb{G}$ and some well known
results that will be used in the sequel. In Section 3 we prove
$L^p-$ Hardy inequality on polarizable Carnot groups. Finally, in
Section 4 we study problem (1.1).

\section {Carnot group}

 A Carnot group ( see \cite {2}, \cite{3},
\cite{13}, \cite{14}, \cite{15}, \cite{31} and \cite{34}) is a
connected, simply connected, nilpotent Lie group
$\mathbb{G}\equiv(\mathbb{R}^n,\cdot)$ whose Lie algebra
$\mathcal{G}$ admits a stratification. That is, there exist linear
subspaces $V_1, ..., V_k$  of $\mathcal{G}$ such that
\begin{equation} \mathcal{G}=V_1\oplus...\oplus V_k, \quad [V_1, V_i]=V_{i+1},
\quad \text{for}\quad i=1,2, ..., k-1 \quad\text {and}\quad [V_1,
V_k]=0
\end{equation}
 where $[V_1, V_i]$ is the subspace  of $\mathcal{G}$
generated by the elements $[X,Y]$ with $X\in V_1$ and $Y\in V_i$.
This defines a $k$-step Carnot group and integer $k\ge 1$ which is
called the step of $ \mathbb{G}$.

Via the exponential map, it is possible to induce on $\mathbb{G}$
a family of automorphisms of the group, called dilations,
$\delta_{\lambda}: \mathbb{R}^n\longrightarrow \mathbb{R}^n
(\lambda>0)$ such that
\[\delta_{\lambda}(x_1, ...,x_n)=(\lambda^{\alpha_1} x_1,...,\lambda^{\alpha_n}x_n)\]
where $1=\alpha_1=...=\alpha_m<\alpha_{m+1}\le ...\le \alpha_n$
are integers and $m=\text{dim}(V_1)$.

The group law can be written in the following form
\begin{equation}
x\cdot y=x+y+P(x,y), \quad x, y\in \mathbb{R}^n
\end{equation}
where $ P:\mathbb{R}^n\times \mathbb{R}^n\longrightarrow
\mathbb{R}^n$ has polynomial components and $P_1=...=P_m=0$. Note
that the inverse $x^{-1}$ of an element $x\in \mathbb{G}$ has the
form $x^{-1}=-x=(-x_1,....-x_n)$.

Let $X_1, . . ., X_m$ be a family of left invariant vector fields
that is an orthonormal basis of $ V_1\equiv\mathbb{R}^m$ at the
origin, that is, $X_1(0)=\partial_{x_1}, . . ., X_m(0)=\partial
_{x_m}$. The vector fields $X_j$  have polynomial coefficients and
can be assumed to be of the form
\[X_j(x)=\partial_j+\sum_{i=j+1}^n a_{ij}(x)\partial_i, \quad
X_j(0)=\partial_j, j=1, . . ., m,\] where  each polynomial
$a_{ij}$ is homogeneous with respect to the dilations of the
group, that is
$a_{ij}(\delta_{\lambda}(x))=\lambda^{\alpha_i-\alpha_j}a_{ij}(x)$.
The horizontal  gradient on $\mathbb{G}$ is the vector valued
operator
\[\nabla_{\mathbb{G}}=(X_1, . . ., X_m)\] where $X_1, . . ., X_m$ are the generators of $ \mathbb{G}$.

The nonlinear operator
\[\Delta_{\mathbb{G},p}\,u=\nabla_{ \mathbb{G}}\cdot (|\nabla_{
\mathbb{G}}u|^{p-2}\nabla_{\mathbb{G}}u)\] is the sub-p-Laplacian
on Carnot group $ \mathbb{G}$. If $p=2$ then we have sub-Laplacian
\[\Delta_{ \mathbb{G}}=\sum_{j=1}^m X_j^2\]
which is a second-order partial differential operator on $
\mathbb{G}$. The fundamental solution $u$ for
$\Delta_{\mathbb{G}}$ is defined to be a weak solution to the
equation
\[-\Delta_{\mathbb{G}}u=\delta\] where $\delta$ denotes the Dirac distribution with singularity at
the neutral element $0$ of $ \mathbb{G}$. In \cite{13},  Folland
proved that in any Carnot group $ \mathbb{G}$,  there exists a
homogeneous norm $N$ such that
\[u=N^{2-Q}\] is a fundamental solution for $\Delta_{
\mathbb{G}}$ ( see also \cite{6}).

We now set
\begin{equation}
N(x):=
\begin{cases}
u^{\frac{1}{2-Q}}&\quad\text{if}\quad x\neq 0,\\
 0 &\quad\text{if}\quad x=0.
 \end{cases}
 \end{equation}  We
recall that a homogeneous norm on $ \mathbb{G}$ is a continuous
 function $N:\mathbb{G}\longrightarrow [0, \infty )$ smooth away from the origin
which satisfies the conditions : $N(\delta_{\lambda}(x))=\lambda
N(x)$, $N(x^{-1})=N(x)$ and $N(x)=0$ iff $x=0$.

 The curve $\gamma:[a,b]\subset \mathbb{R}\longrightarrow
\mathbb{G}$ is called horizontal if its tangents lie in $V_1$,
i.e,  $\gamma'(t)\in \text{\textit{span}}\{X_1, . . ., X_m\}$ for
all $t$. Then, the Carnot-Car\'ethedory distance $d_{CC}(x,y)$
between two points $x,y\in \mathbb{G}$ is defined to be the
infimum of all horizontal lengths $\int_a^b \langle \gamma'(t),
\gamma'(t)\rangle^{1/2} dt$ over all horizontal curves
$\gamma:[a,b]\longrightarrow \mathbb{G}$ such that $\gamma(a)=x$
and $\gamma(y)=b$. Notice that $d_{CC}$ is a homogeneous norm and
satisfies the invariance property

\[d_{CC}(z\cdot x, z\cdot y)=d_{CC}(x,y), \quad \forall\, x,y,
z\in \mathbb{G},\] and is homogeneous of degree one with respect
to the dilation $\delta_{\lambda}$, i.e.
\[d_{CC}(\delta_{\lambda}(x), \delta_{\lambda}(y))=\lambda
d_{CC}(x,y), \quad \forall\, x, y, z \in \mathbb{G}, \forall\,
\lambda>0.\]

The Carnot-Careth\'edory balls are defined by $B(x, R)=\{y\in
\mathbb{G} |d_{CC}(x,y)<R\}$.  The $n$-dimensional Lebesgue
measure  $\mathcal{L}^n$, is the Haar measure of group $
\mathbb{G}$. This means that if $E \subset \mathbb{R}^n$ is a
measurable, then $\mathcal{L}^n(x\cdot E)=\mathcal{L}^n(E)$ for
all $x\in \mathbb{G}$. Moreover, if $\lambda>0$ then $
\mathcal{L}^n(\delta_{\lambda}(E))=\lambda^Q \mathcal{L}^n(E)$.
Clearly
\[ \mathcal{L}^n\big( B(x, R)\big)=R^Q \mathcal{L}^n\big( B(x,1)\big)=R^Q \mathcal{L}^n\big(B(0,1)\big)\] where \[Q=\sum_{j=1}^k j (\text{dim}V_j)\] is the
homogeneous dimension of $ \mathbb{G}$.

It is well known that  Sobolev inequalities are important in the
study of partial differential equations, especially in the study
of those arising from geometry and physics. The following Sobolev
inequality holds  on $ \mathbb{G}$
\begin{equation}
\Big(\int_{\mathbb{G}}|\phi(x)|^qdx\Big)^{1/q}\le
C_{p,q}\Big(\int_{\mathbb{G}}|\nabla_{\mathbb{G}}\phi(x)|^pdx\Big)^{1/p}
\end{equation}
where $1\le p<Q$, $q=\frac{Qp}{Q-p}$ ( see \cite {13}, \cite{35} )
and plays a crucial role in this paper.
\medskip

\noindent \textbf{Example 1.} The simplest non trivial example of
Carnot group is given by the Heisenberg group $ \mathbb{H}^n$.
Denoting points in $ \mathbb{H}^n$ by $(z,t)$ with $z=(z_1, . . .,
z_n)\in \mathbb{C}^n$ and $t\in \mathbb{R}$ we have the group law
given as
\[(z,t)\circ (z',t')=(z+z', t+t'+2\sum_{j=1}^n
Im(z_j\bar{z}_j'))\] With the notation $z_j=x_j+iy_j$, the
horizontal space $V_1$ is spanned by the basis
\[X_j=\frac{\partial}{\partial x_j}
+2y_j\frac{\partial }{\partial t}\quad\text{and} \quad
Y_j=\frac{\partial}{\partial y_j}-2x_j\frac{\partial }{\partial
t}.\] The one dimensional center $V_2$ is spanned by the vector
field $T=\frac{\partial}{\partial t}$.  We have the commutator
relations $[X_j, Y_j]=-4T$, and all other brackets of $\{X_1, Y_1,
. . ., X_n, Y_n\}$ are zero. The sub-Laplacian associated  with
the basis $\{X_1, Y_1, ..., X_n, Y_n\}$  is the operator
\[\Delta_{ \mathbb{H}^n}=\sum_{j=1}^n (X_j^2+Y_j^2).\] A
homogeneous norm on $ \mathbb{H}^n$ is given by
\[\rho=|(z,t)|=(|z|^4+t^2)^{1/4}\] and the homogeneous dimension of $\mathbb{H}^n$ is $Q=2n+2$.

A remarkable analogy between sub-Laplacian and the classical
Laplace operator has been obtained by Folland \cite{12}. He found
that the fundamental solution of $-\Delta_{ \mathbb{H}^n}$ with
pole zero is given by
\[\Psi(z,t)=\frac{c_Q}{\rho(z,t)^{Q-2}}\quad\text{where}\quad
c_Q=\frac{2^{(Q-2)/2}\Gamma((Q-2)/4)^2}{\pi^{Q/2}}.\]

\medskip

\noindent \textbf{Example 2.} Another important model of Carnot
groups are the \textit{H}-type (Heisenberg type) groups which were
introduced by Kaplan \cite {25} as direct generalizations of the
Heisenberg group $\mathbb{H}^n$. An \textit{H}-type group is a
Carnot group with a two-step Lie algebra $ \mathcal{G}=V_1\oplus
V_2$ and an inner product $\langle, \rangle$ in $ \mathcal{G}$
such that the linear map
\[ J:V_2\longrightarrow \text{End}V_1,\] defined by the condition
\[\langle J_z(u), v\rangle=\langle z, [u,v]\rangle,\quad u,v\in V_1, z\in V_2\] satisfies
\[J_z^2=-||z||^2\mathbf{Id}\]
for all $z\in V_2$, where $||z||^2= \langle z,z\rangle$.

Sub-Laplacian is defined in terms of a fixed basis $X_1,. . . ,
X_m$ for $V_1$:

\begin{equation}
\Delta_\mathbb{G}=\sum_{i=1}^mX_i^2.
\end{equation}

The exponential mapping of a simply connected Lie group is an
analytic diffeomorphism. One can then define analytic mappings
$v:\mathbb{G} \longrightarrow V_1$ and
$z:\mathbb{G}\longrightarrow V_2$ by
\[x=\text{exp}(v(x)+z(x))\] for every $x\in \mathbb{G}$.
In \cite {25}, Kaplan proved that there exists a constant $c>0$
such that the function
\[\Phi (x)=c\Big(|v(x)|^4+16|z(x)|^2\Big)^{\frac{2-Q}{4}}\]
is a  fundamental solution for the operator $ -\Delta_\mathbb{G}$.
We note that
\begin{equation}
K(x)=\Big(|v(x)|^4+16|z(x)|^2\Big)^{\frac{1}{4}}
\end{equation}
defines a homogeneous norm  and $Q=m+2k$ is the homogeneous
dimension of $\mathbb{G}$ where $m=$dim$V_1$ and $k=$dim$V_2$.
This result generalized Folland's fundamental solution for the
Heisenberg group $ \mathbb{H}^n$ \cite{12}. It is remarkable that
the homogeneous norm $K(x)$ (2.6) involves also in the expression
of fundamental solution of the following $p-$sub-Laplace operator
\begin{equation}
\mathcal{L}_pu=\sum_{i=1}^mX_i(|Xu|^{p-2}X_iu), \quad 1<p<\infty.
\end{equation}
More precisely, Capogna, Danielli and Garofalo  proved that for
every $1<p<\infty$ there exists $c_p>0$ and $\tilde{c}_p<0$ such
that the function
\begin{equation}
\Gamma_p(x)=
\begin{cases}
c_pK^{(p-Q)/(p-1)} &\quad \text{when}\, p\ne Q,\\
\tilde{c}_p\log K&\quad\text{when}\, p=Q,
\end{cases}
\end{equation}
is a fundamental solution for the operator $-\mathcal{L}_p$
\cite{9}. This nice result immediately motivate the following
question: Is it possible to find explicit formulas for the
fundamental solutions to the sub-$p-$Laplacian on more general
Carnot groups? An interesting  result in this direction has been
obtained by Balogh and Tyson \cite{2}. They found explicit
formulas for the fundamental solutions to the sub-$p-$Laplacian on
polarizable Carnot goups.

\section{Poloraziable Carnot groups and Hardy type inequality}
\noindent \textbf{Polarizable Carnot group.} A Carnot
 group $ \mathbb{G}$ is called polarizable if the homogeneous norm $N=u^{1/(2-Q)}$,
 associated to Folland's solution $u$ for the sub-Laplacian
 $\Delta_{\mathbb{G}}$, satisfies the following $\infty$- sub-Laplace
 equation,

 \begin{equation}
\Delta_{\mathbb{G},\infty}N:= \frac{1}{2}\langle\nabla_{
\mathbb{G}}(|\nabla_{ \mathbb{G}}N|^2), \nabla_{
\mathbb{G}}N\rangle=0, \quad\quad \text{in} \quad
\mathbb{G}\setminus \{0\}.
 \end{equation}
In \cite {2},  Balogh and Tyson proved that the fundamental
solutions of the sub-$p$-Laplacian
\begin{equation}
\Delta_{ \mathbb{G}, p}f=\text{div}(|\nabla_{
\mathbb{G}}f|^{p-2}\nabla_{\mathbb{G}}f), \quad\quad f\in C^2(
\mathbb{G}),
\end{equation}
 is given
by
\begin{equation}
u_p=
\begin{cases}
\begin{aligned}  N^{\frac{p-Q}{p-1}}, \quad\text{if}\quad p \neq
Q\\
-\text{log}N, \quad\text{if}\quad p=Q.
\end{aligned}
\end{cases}
\end{equation}
Moreover, for each $1<p<\infty$ there exists a constants $c_p>0$
so that the function $c_pu_p$ solves the equation $-\Delta_p
(c_pu_p)=\delta$ in the sense of distributions. This class of
groups is the largest class for which the fundamental solution of
the sub-$p$-Laplacian (3.2) is given by (3.3) \cite {2}.

It is known that Euclidean space, Heisenberg group and the
Kaplan's H-type group are polarizable Carnot group. Further
information about polarizable Carnot groups can be found in
\cite{2}.
\medskip

\noindent \textbf{Hardy-type inequalities.} To motivate our
results, let us recall the following  $L^p-$ Hardy inequality in
Euclidean space $\mathbb{R}^n$:
\begin{equation} \int_{\mathbb{R}^n} |\nabla\phi(x)|^p dx\ge
|\frac{n-p}{p}|^p \int_{\mathbb{R}^n}
\frac{|\phi(x)|^p}{|x|^p}dx\end{equation}  where  $\phi\in
C_0^{\infty}( \mathbb{R}^n\setminus \{0\})$ and
$|\frac{n-p}{p}|^p$ is the best constant. The inequality was
discovered by Hardy \cite {24} in the one-dimensional case. Later
on it has been extended to higher dimensions and generalized
involving various kinds of distance functions in Euclidean space
(see \cite {5}, \cite{32} and the references therein)

In recent years there has been considerable progress in the study
of Hardy-type inequalities on sub-Rimeannian spaces. The following
inequality is the analogue of the inequality (3.4) on the
Heisenberg group $ \mathbb{H}^n$

\begin{equation}
\int_{\mathbb{H}^n} |\nabla_{ \mathbb{H}^n}\phi(z, l)|^pdzdl\ge
\big(\frac{Q-p}{p}\big)^p \int_{\mathbb{H}^n}
\frac{|z|^p}{(|z|^4+l^2)^{\frac{p}{2}}}|\phi(z,l)|^pdzdl
\end{equation}
where $Q=2n+2$, $1<p<Q$ and $\phi\in C_0^{\infty}(
\mathbb{R}^n\setminus \{0\})$.  The inequality (3.5) was first
proved by Garofalo and Lanconelli \cite {17} for $p=2$ and later
by Niu, Zhang and Wang \cite {30} for any $p$. ( See also
\cite{10}, \cite{22}).
\medskip

These inequalities play an important role in the study of linear
and nonlinear partial differential equations on Euclidean space
$\mathbb{R}^n$ and on the Heisenberg group $ \mathbb{H}^n$ (see
\cite{1}, \cite{4}, \cite{7}, \cite{8}, \cite{11}, \cite{16},
\cite {19}, \cite{20}, \cite{21}, \cite{22}, \cite{23}, \cite{26},
\cite{27}, \cite{33}).
\medskip

In this paper we prove $L^p-$ Hardy inequality  on polarizable
Carnot groups and it is an open problem for more general Carnot
groups. Recently, Kombe has proved Hardy inequality with a sharp
constant for $p=2$ on general Carnot groups \cite{28}. The
following theorem is the main result of this section.
\begin{theorem} Let $\mathbb{G}$ be a polarizable Carnot group with homogeneous norm  $N=u^{1/(2-Q)}$  and  let $\phi\in
C_0^{\infty}(\mathbb{G}\setminus \{0\})$, $ Q\ge 3$ and $1<p<Q$.
Then the following inequality is valid
\begin{equation}\int_{\mathbb{G}} |\nabla_{\mathbb{G}} \phi|^p dx
\ge (\frac{Q-p}{p})^p \int_{\mathbb{G}} \frac{|\nabla_{\mathbb{G}}
N|^p}{N^p}|\phi| ^pdx.\end{equation} Furthermore,  the constant
$(\frac{Q-p}{p})^p$ is sharp.
\end{theorem}
\proof

Let $\phi=N^{\gamma}\psi$  where $\psi\in
C_0^{\infty}(\mathbb{G}\setminus \{0\})$ and $\gamma\in
\mathbb{R}\setminus \{0\}$. A direct calculation
 shows that
 \begin{equation}
|\nabla_{\mathbb{G}}(N^{\gamma}\psi)|=|\gamma
N^{\gamma-1}\psi\nabla_{\mathbb{G}} N+
N^{\gamma}\nabla_{\mathbb{G}}\psi|.
\end{equation}
We now use the following inequality which is  valid for any $a, b
\in \mathbb{R}^n$ and  $1<p<2$,
\begin{equation}
|a+b|^p-|a|^p\ge
c(p)\frac{|b|^2}{(|a|+|b|)^{2-p}}+p|a|^{p-2}a\cdot b
\end{equation}
where $c(p)> 0$ (see \cite{5}, \cite {29} ).
 In view of (3.8) we have
that
\[\begin{aligned} \int_{\mathbb{G}} |\nabla_{\mathbb{G}} \phi|^p dx
 &\ge |\gamma|^p \int_{\mathbb{G}}N^{\gamma
p-p}|\nabla_{\mathbb{G}}N|^p|\psi|^pdx+|\gamma|^{p-2}\gamma
\int_{\mathbb{G}} N^{\gamma p-p+1}|\nabla_{\mathbb{G}}
N|^{p-2}\nabla_{\mathbb{G}} N \cdot
\nabla_{\mathbb{G}}(|\psi|^p)dx \\ &+ c(p)\int_{\mathbb{G}}\frac{
N^{2\gamma}|\nabla_{\mathbb{G}}\psi|^2}{(|\gamma
N^{\gamma-1}\psi\nabla_{\mathbb{G}} N|+
|N^{\gamma}\nabla_{\mathbb{G}}\psi|)^{2-p} }dx.\end{aligned}\]
Clearly
\[
\int_{\mathbb{G}} |\nabla_{\mathbb{G}} \phi|^p dx
 \ge |\gamma|^p \int_{\mathbb{G}}N^{\gamma
p-p}|\nabla_{\mathbb{G}}N|^p|\psi|^pdx +|\gamma|^{p-2}\gamma
\int_{\mathbb{G}} N^{\gamma p-p+1}|\nabla_{\mathbb{G}}
N|^{p-2}\nabla_{\mathbb{G}} N \cdot
\nabla_{\mathbb{G}}(|\psi|^p)dx,
\] and the integration by parts gives
\[
\int_{\mathbb{G}} |\nabla_{\mathbb{G}} \phi|^p dx
 \ge |\gamma|^p \int_{\mathbb{G}}N^{\gamma
p-p}|\nabla_{\mathbb{G}}N|^p|\psi|^pdx-|\gamma|^{p-2}\gamma
\int_{\mathbb{G}} \nabla_{\mathbb{G}}\cdot( N^{\gamma
p-p+1}|\nabla_{\mathbb{G}} N|^{p-2}\nabla_{\mathbb{G}} N)
|\psi|^pdx.
\]
We now choose $\gamma=\frac{p-Q}{p}$ then we get
\begin{equation}
\begin{aligned} \int_{\mathbb{G}} \nabla_{\mathbb{G}}\cdot( N^{\gamma p-p+1
}|\nabla_{\mathbb{G}} N|^{p-2}\nabla_{\mathbb{G}} N) |\psi|^pdx
&=\int_{\mathbb{G}} \nabla_{\mathbb{G}}\cdot(
N^{1-Q}|\nabla_{\mathbb{G}} N|^{p-2}\nabla_{\mathbb{G}} N)
|\psi|^pdx\\
&=\int_{\mathbb{G}} (\Delta_{\mathbb{G},p}u_p)|\psi|^pdx.
\end{aligned}
\end{equation}
Since $u_p$ is the fundamental solution of sub-p-Laplacian
$-\Delta_{\mathbb{G},p}$, then we have
\begin{equation}
\begin{aligned} \int_{\mathbb{G}} \nabla_{\mathbb{G}}\cdot(
N^{1-Q}|\nabla_{\mathbb{G}} N|^{p-2}\nabla_{\mathbb{G}} N)
|\psi|^pdx
&=\int_{\mathbb{G}} (\Delta_{\mathbb{G},p}u_p)|\psi|^pdx\\
&=-|\phi(0)|^pN^{(Q-p)}(0)\\&=0.
\end{aligned}
\end{equation}
Therefore,
\begin{equation}
\int_{\mathbb{G}} |\nabla_{\mathbb{G}} \phi|^p dx
 \ge \Big(\frac{Q-p}{p}\Big)^p
 \int_{\mathbb{G}}\frac{|\nabla_{\mathbb{G}} N|^p}{N^p}|\phi|^pdx.
\end{equation}
The theorem (3.1) also holds for $p> 2$ and in this case we use
the following inequality
\begin{equation}
|a+b|^p-|a|^p\ge c(p)|b|^p+p|a|^{p-2}a\cdot b
\end{equation}
where $a\in  \mathbb{R}^n, b\in \mathbb{R}^n$ and $ c(p)> 0$ (see
\cite{5}, \cite{29} ). If $p=2$ then the inequality (3.11) holds
in any Carnot group without the hypothesis of polarizability \cite
{28}.
\medskip

To show that the constant $\Big(\frac{Q-p}{p}\Big)^p$ is sharp, we
use the following family of functions
\begin{equation} \phi_{\epsilon}(x)=
\begin{cases}
 1 &\quad\text{if}  \quad N(x)\in [0,1],\\
N^{-(\frac{Q-p}{p}+\epsilon)} &\quad \text{if} \quad N(x)>1,
\end{cases}
\end{equation}
and pass to the limit as $\epsilon\longrightarrow 0$. Here we
notice that $|\nabla_{\mathbb{G}}N|$ is uniformly bounded  and
polar coordinate integration formula holds
 on  $\mathbb{G}$ \cite {2}, \cite{14}.
\section{Nonexistence results}
We now turn our attention to problem (1.1). Throughout this
section, $\Omega$ is a bounded domain in  $ \mathbb{G}$ with
smooth boundary and $V\in L_{\text{loc}}^1(\Omega\setminus
\mathcal{K})$ for some closed Lebesque null set $ \mathcal{K}$. We
define the positive solutions in the following sense.

\begin{definition} By a $\textit{positive local solution
continuous}$ off of $\mathcal{K}$, we mean
\begin{itemize}
  \item [(i)] $ \mathcal{K}$ is a closed Lebesgue null subset of $\Omega$,
  \item [(ii)] $u :[0, T)\longrightarrow L^1(\Omega)$ is
  continuous for some $T>0$,
  \item [(iii)] $(x,t)\longrightarrow u(x,t)\in C((\Omega\setminus
  \mathcal{K})\times (0,T))$,
 \item [(iv)] $u(x,t)>0$ on $(\Omega\setminus \mathcal{K})\times (0,T)$,
 \item [(v)] $\lim_{t \to 0}u(.,t)=u_0 $ in the sense of
 distributions,
 \item [(vi)]$\nabla_{ \mathbb{G}} u \in L_{loc}^p(\Omega\setminus \mathcal{K})$  and $u$ is a solution in the
 sense of distributions  of the PDE.
\end{itemize}
\noindent {\bf Remark.} If  $0<a<b<T$ and $\mathcal{K}_o$ is a
compact subset of $\Omega\setminus  \mathcal{K}$, then
$u(x,t)\ge\epsilon_1>0$ for $(x,t)\in \mathcal{K}_o\times [a,b]$
for some $\epsilon_1>0$. We can weaken (iii), (iv) to be
\begin{itemize}
  \item [(iii)'] $u(x,t)$ is positive and locally bounded on
  $(\Omega\setminus \mathcal{K})\times (0,T)$,
  \item [(iv)'] $\frac{1}{u(x, t)}$ is locally bounded on
  $(\Omega\setminus \mathcal{K})\times (0,T)$.
\end{itemize}
\end{definition}
If a solution satisfies  (i), (ii), (iii)', (iv)', (v), and (vi)
then  we call it a `` general positive local solution off of
$\mathcal{K}$ ". This is more general than a positive local
solution continuous off of $ \mathcal{K}$. If $ \mathcal{K}=
\emptyset$, we simply call $u$ ``general positive local solution".
\begin{theorem}
Let $\frac{2Q}{Q+1}\le p<2 $ and $V\in
L_{\text{loc}}^1(\Omega\setminus \mathcal{K})$ where $
\mathcal{K}$ is a Lebesgue null subset of $\Omega$. If
\begin{equation} \sigma_{\inf}^p((1-\epsilon)V) :=\inf_{0 \neq \phi \in
C_0^{\infty}(\Omega\setminus \mathcal{K})}
\frac{\int_{\Omega}|\nabla_{\mathbb{G}}\phi|^pdx
-\int_{\Omega}(1-\epsilon) V |\phi|^p dx}{\int_{\Omega}|\phi|^p
dx}=-\infty \end{equation} for some $\epsilon >0 $, then (1.1) has
no general positive local solution off of $ \mathcal{K}$.
\end{theorem}
\proof  We argue by contradiction. Given any $T>0$, let $u: [0,
T)\longrightarrow L^1(\Omega)$ be a general positive local
solution to (1.1) in $(\Omega\setminus \mathcal{K})\times (0, T)$
with $u_0 \ge 0$ but not identically zero. Multiply both sides of
(1.1) by the test function $|\phi|^p/u^{p-1}$ where $\phi\in
C_0^{\infty}(\Omega\setminus \mathcal{K})$, and integrate over
$\Omega$, to get

\begin{equation}
\frac{1}{2-p}\frac{\partial }{\partial t}\int_{\Omega}u^{2-p}
|\phi|^p dx-\int_{\Omega}\nabla_{
\mathbb{G}}\cdot(|\nabla_{\mathbb{G}}u|^{p-2}\nabla_{\mathbb{G}}
u) (\frac{|\phi|^p}{u^{p-1}})dx = \int_{\Omega}V|\phi|^p dx.
\end{equation}
It follows from the integration by parts that
\begin{equation}
L=\int_{\Omega}\nabla_{\mathbb{G}}\cdot(|\nabla_{\mathbb{G}}u|^{p-2}\nabla_{\mathbb{G}}
u) (\frac{|\phi|^p}{u^{p-1}})dx
=-\int_{\Omega}|\nabla_{\mathbb{G}} u|^{p-2} \nabla_{\mathbb{G}}
u\cdot \nabla_{\mathbb{G}}(\frac{|\phi|^p}{u^{p-1}})dx.
\end{equation}
Since
\[|\nabla _{\mathbb{G}}u|^{p-2}\nabla_{
\mathbb{G}} u\cdot\nabla_{
\mathbb{G}}(\frac{|\phi|^p}{u^{p-1}})=p|\nabla_{\mathbb{G}}
u|^{p-2} \frac{|\phi|^{p-1}}{u^{p-1}}\nabla_{
\mathbb{G}}u\cdot\nabla_{\mathbb{G}}
|\phi|-(p-1)\frac{|\phi|^p}{u^p}|\nabla _{\mathbb{G}}u|^p,
\]
\begin{equation}
L = (p-1)\int_{\Omega}|\nabla _{\mathbb{G}} u|^p
\frac{\phi^p}{u^p}dx-\text{sign}_0(\phi)p\int_{\Omega}|\nabla_{\mathbb{G}}u|^{p-2}
\frac{\phi^{p-1}}{u^{p-1}}\nabla_{\mathbb{G}}u\cdot\nabla_{\mathbb{G}}|\phi|
dx,
\end{equation}
where \[ \text{sign}_0(\phi)=
\begin{cases}
\frac{\phi}{|\phi|} &\text{if}\quad \phi\neq 0, \\
0 &\text{if}\quad \phi= 0.
\end{cases}\]
Therefore,

\begin{equation}
L \ge (p-1)\int_{\Omega}|\nabla_{ \mathbb{G}}u|^p
\frac{|\phi|^p}{u^p}dx-p \int_{\Omega}|\nabla _{\mathbb{G}}
u|^{p-1}|\nabla_{\mathbb{G}} \phi| \frac{\phi^{p-1}}{u^{p-1}} dx.
\end{equation}
Here we can use the following elementary inequality: Let $p>1$ and
$w_1\neq w_2$ be two positive  real numbers. Then

\[w_1^p-w_2^p-pw_2^{p-1}(w_1-w_2)> 0; \]it follows that

\[(p-1)w_2^p-pw_2^{p-1}w_1> -w_1^p.\]
We can take $w_2=|\frac{\phi}{u}\nabla_{\mathbb{G}} u|$,
$w_1=|\nabla_{\mathbb{G}} \phi|$; then  we have
\begin{equation}
(p-1)\int_{\Omega}|\nabla_{ \mathbb{G}} u|^p
\frac{|\phi|^p}{u^p}dx -p \int_{\Omega}|\nabla_{\mathbb{G}}
u|^{p-1}|\nabla_{\mathbb{G}} \phi| \frac{|\phi|^{p-1}}{u^{p-1}}
dx\ge -\int_{\Omega}|\nabla_{\mathbb{G}} \phi|^pdx.
\end{equation}
Therefore

\begin{equation}
L=\int_{\Omega}\nabla_{\mathbb{G}}\cdot(|\nabla_{\mathbb{G}}u|^{p-2}\nabla_{\mathbb{G}}
u) (\frac{|\phi|^p}{u^{p-1}})dx \ge -\int_{\Omega}|\nabla_{
\mathbb{G}}\phi|^p dx.
\end{equation}
Substituting  (4.7) into (4.2) and integrating from $t_1$ to
$t_2$, where $0<t_1<t_2<T$, we obtain
\begin{equation}
\int_{\Omega}V(z)|\phi|^p dx-\int_{\Omega}|\nabla_{\mathbb{G}}
\phi|^p dx \le
\frac{1}{(2-p)(t_2-t_1)}\int_{\Omega}(u^{2-p}(x,t_2)-u^{2-p}(x,t_1))|\phi|^p
dx .
\end{equation} Using Jensen's inequality for concave
functions, we obtain
\[\int_{\Omega}\Big (u(x, t_i)\Big )^{\frac{(2-p)(Q)}{p}}dx \le
C(|\Omega|) \Big(\int_{\Omega} u(x,t_i) dx \Big
)^{\frac{(2-p)(Q)}{p}}<\infty.\] Therefore

\[u^{2-p}(x,t_i) \in L^{\frac{Q}{p}}(\Omega).\]
We now use the following a priori inequality which is a
consequence of the Sobolev inequality (2.4). For every
$\epsilon>0$ there exists $C(\epsilon)$ such that
\begin{equation}
\begin{aligned}
&\frac{1}{(2-p)(t_2-t_1)}\int_{\Omega}\Big(u^{2-p}(x,t_2)-u^{2-p}(x,t_1)\Big)|\phi|^p
dx \\ &\le
\frac{\epsilon}{1-\epsilon}\int_{\Omega}|\nabla_{\mathbb{G}}
\phi|^pdx+C(\epsilon)\int_{\Omega}|\phi|^2dx.
\end{aligned}
\end{equation}
Substituting (4.9) into (4.8), we obtain
\begin{equation}
\int_{\Omega}V(z)|\phi|^p dx-\int_{\Omega}|\nabla_{\mathbb{G}}
\phi|^p dx\le
\frac{\epsilon}{1-\epsilon}\int_{\Omega}|\nabla_{\mathbb{G}}\phi|^p
dx +C(\epsilon)\int_{\Omega}|\phi|^pdx.
\end{equation}
Therefore
\begin{equation}
\inf_{0 \neq \phi \in C_c^{\infty} (\Omega\setminus
\mathcal{K})}\frac{ \int_{\Omega}|\nabla_{ \mathbb{G}}\phi|^p
dx-\int_{\Omega}(1-\epsilon)V(x)|\phi|^p dx
}{\int_{\Omega}|\phi|^p dx} \ge -(1-\epsilon)C(\epsilon)>-\infty.
\end{equation}
This contradicts our assumption (4.1). The proof of Theorem 4.1 is
now complete. $\hfill\qed$
\medskip

\noindent\textbf{Singular Potentials.}  We now focus our attention
on some singular potentials. First we treat the following positive
singular potential: \begin{equation} V(x)=\lambda\frac{|\nabla_{
\mathbb{G}}N|^p}{N^p}\end{equation} where $\lambda>0$.

As a sign changing potential, we consider the following highly
singular, oscillating potential:
\begin{equation} V(x)=\lambda\frac{|\nabla_{ \mathbb{G}}N|^p}{N^p}+\beta
\frac{|\nabla_{ \mathbb{G}}N|^p}{N^p}\sin
(\frac{1}{N^{\alpha}})\end{equation}  where $ \lambda>0$,
$\beta\in \mathbb{R}\setminus\{0\}$ and $\alpha\in \mathbb{R}$.

Theorem 4.1 gives the following corollaries.
\begin{corollary}
Let $0\in\Omega$ and $V(x)$ be defined by (4.12). Then the problem
(1.1) has no general positive local solution off of $ \mathcal{K}$
if $\frac{2Q}{Q+1}\le p<2$ and $\lambda>(\frac{Q-p}{p})^p$.
\end{corollary}
\begin{corollary}
Let $0\in\Omega$ and $V(x)$ be defined by (4.13). Then the problem
(1.1) has no general positive local solution off of $ \mathcal{K}$
if $\frac{2Q}{Q+1}\le p<2$ and $\lambda>(\frac{Q-p}{p})^p$.
\end{corollary}

To prove  Corollary 4.1 and Corollary 4.2, we need to construct a
sequence $\{\phi_n\}$ of test functions which satisfies
\begin{equation}
\frac{\int_{\Omega}|\nabla_{ \mathbb{G}} \phi_n|
^pdx-\int_{\Omega}(1-\epsilon) V(x) |\phi_n|^p(x)
dx}{\int_{\Omega}|\phi_n|^pdx}\longrightarrow
-\infty\quad\text{as}\quad n\longrightarrow\infty.
\end{equation}
We can easily construct
 radial functions $\phi_n=f(N(x))$ which satisfy 4.13. (see \cite
 {20}, \cite{21}).
\medskip

Notice that when $|\beta|>\lambda$, the potential in Corollary 4.2
has very large positive and negative parts, in particular, it
oscillates wildly, but important cancellations occur between the
positive and the negative parts in the quadratic form. Therefore,
nonexistence of positive solutions only depends on the size of
$\lambda$.
\medskip

\noindent {\bf Remark.} We would like to point out that all the
new results in Section 4 hold on the whole space $\mathbb{G}$
which is a consequence of bounded domain case (see \cite{23} or
\cite{26}).
\bibliographystyle{amsalpha}

\end{document}